\documentclass[12pt]{amsart}%
\usepackage{amsmath,amsthm,amsfonts,amssymb}
\usepackage{verbatim}
\usepackage{hyperref}
\usepackage{latexsym}
\usepackage{marginnote}
\usepackage{amscd}
\usepackage{amssymb}
\usepackage{xcolor}
\usepackage{amsmath}
\usepackage{amsfonts}
\usepackage{mathrsfs}
\usepackage{graphicx}%
\providecommand{\U}[1]{\protect\rule{.1in}{.1in}}
\providecommand{\U}[1]{\protect\rule{.1in}{.1in}}
\providecommand{\U}[1]{\protect\rule{.1in}{.1in}}
\newtheorem{theorem}{Theorem}[section]

\newtheorem{definition}[theorem]{Definition}

\newtheorem{example}[theorem]{Example}

\newtheorem{remark}[theorem]{Remark}

\def\Var{\mathop{Var}}
\begin{document}
\title[On Spectral Barron spaces on compact groups ]{Spectral Barron spaces of vector-valued functions on compact groups}
\author{Yaogan Mensah and Isiaka Aremua}

\address{Department of Mathematics, University of Lom\'e, Togo}
\email{\textcolor[rgb]{0.00,0.00,0.84}{mensahyaogan2@gmail.com}}

\address{Department of Physics, University of Lom\'e, Togo}
\email{\textcolor[rgb]{0.00,0.00,0.84}{claudisak@gmail.com}}

\begin{abstract}
 In this article, we study spectral Barron spaces whose elements are
Banach space-valued functions on a compact group whose Fourier transforms
admit a certain summability property. We investigate the functional
properties of these spaces and establish continuous embeddings with respect
to other function spaces, among which are Sobolev spaces of vector-valued
functions and the space of bounded vecto-valued functions on compact
groups. Beyond these structural results, we prove a quantitative
approximation theorem: when the target space is a separable Hilbert space,
every function in the spectral Barron space  admits an approximation by matrix coefficients of unitary representations of the group with $L^2$-error decaying at the rate $O(n^{-1/2})$ and a constant
depending only on the spectral Barron norm of the function. This extends,
via Maurey's empirical method, the classical dimension-independent
approximation rate of Barron's theorem beyond the Euclidean setting, and we
further indicate how the argument persists, with an inflated constant, when
target space is only assumed to have Rademacher type 2.
\end{abstract}
\maketitle


\section{Introduction}

Spectral Barron spaces are an emerging class of function spaces in
functional analysis and approximation theory. They provide a rich framework
for controlling neural network approximation error, and they originated
from the seminal papers of A.R. Barron \cite{Barron1993,Barron1994}. In
Barron's original setting on $\mathbb{R}^d$, a function whose Fourier
transform is summable against a suitable weight can be approximated by an
$n$-term superposition of sigmoidal ridge functions with $L^2$-error of
order $O(n^{-1/2})$, with a constant depending only on the Barron norm of
the function and not on the ambient dimension $d$. This
dimension-independent rate is the central feature that makes Barron spaces
useful for high-dimensional approximation, and it rests on a general
Hilbert-space argument known as Maurey's empirical method
\cite{Pisier1981}, later popularized in the approximation-theory literature
by Jones \cite{Jones1992} and adapted by Pisier to Banach spaces of
Rademacher type 2 \cite{Pisier1989}.
 
Some recent works on spectral Barron spaces on Euclidean spaces are
\cite{ChoulliLuTakase2025, FengLu2025, LiaoMing2025, SiegelXu2023}. This
list is far from exhaustive. In \cite{LiaoMing2025}, Liao and Ming proved a
sharp embedding between the spectral Barron space and the Besov space with
embedding constants independent of the input dimension. In
\cite{ChoulliLuTakase2025}, Choulli, Lu and Takase conducted a systematic
study of spectral Barron spaces on $\mathbb{R}^d$ with applications to
partial differential equations.
 
The aim of this work is to move beyond the Euclidean framework and to place
spectral Barron spaces at the core of abstract and applied harmonic
analysis. In this article, we study spectral Barron spaces of
Banach space-valued functions on compact groups. To do so, we use the
Fourier transformation of vector-valued functions on compact groups as
formulated by Assiamoua and Olubummo \cite{AssiamouaOlubummo1989}. We first
study the functional properties of these spaces and their continuous
embeddings with respect to Sobolev spaces of vector-valued functions and the
space of bounded vector-valued functions on compact groups. We then go
further than the purely functional-analytic side of the theory:
in Section 4, we establish an explicit, quantitative approximation theorem
for spectral Barron spaces on compact groups. The uniform boundedness of
matrix coefficients of unitary representations allows every $f \in \mathfrak{B}^s_\gamma(G,A)$ to be
written as a barycenter of a bounded dictionary of such matrix coefficients;
Maurey's empirical method then produces, for every $n$, an explicit
$n$-term approximant $f_n$ built from matrix coefficients with
$\|f-f_n\|_{L^2(G,A)} \le \frac{\|f\|_{\mathfrak{B}^s_\gamma}}{\sqrt n}$. This is, to our
knowledge, the first extension of Barron's dimension-independent
approximation rate to functions on general compact groups, with the matrix
coefficients $u^\sigma_{ij}$  playing the role that bounded
sigmoidal ridge functions play on $\mathbb{R}^d$.
 
The rest of the paper is organized as follows. Section \ref{Harmonic analysis of vector-valued function on compact groups} gathers
preliminary notes on Fourier analysis of vector-valued functions on compact
groups and on Assiamoua spaces. Section \ref{Functional and embedding results} contains our main functional and
embedding results for spectral Barron spaces. Section \ref{approximation} contains the Barron-type approximation theorem described above, together with remarks on
the dimension-independence of the rate and its extension to Banach spaces of Rademacher type 2.

\section{Harmonic analysis of vector-valued functions on compact groups}\label{Harmonic analysis of vector-valued function on compact groups}

Let $G$ denote a compact group assumed to be Hausdorff once for all. A unitary representation of the group $G$ on a Hilbert space $H$ is a  homomorphism $\sigma$ from  $G$ into the group $U(H)$ of unitary operators on $H$. The Hilbert space $H$ is called the representation space of $\sigma$ and often denoted $H_\sigma$.  The dimension of $H_\sigma$ is denoted by $d_\sigma$. The representation $\sigma$ is said to be continuous if for each $\xi\in H_\sigma$,  the mapping $G\rightarrow H_\sigma, x\mapsto \sigma(x)\xi$ is continuous. All the representations considered in this paper are assumed to be continuous.

The representation is said to be irreducible if the only closed  subspaces  of $H_\sigma$ that are invariant by $\sigma$ are the zero vector space and the Hilbert space $H_\sigma$ itself. 
It is well known that an irreducible unitary representation of a compact group is necessarily finite-dimensional. 

Two unitary  representations $\sigma_1$ and   $\sigma_2$ are said to be (unitarily) equivalent if there exists a unitary operator $T : H_{\sigma_1} \rightarrow H_{\sigma_2}$ such that 
$$\forall x\in G, \, T\sigma_1 (x)=\sigma_2 (x)T.$$
In the sequel, let $\widehat{G}$ denote the set of equivalent classes of unitary irreducible representations of $G$. The class of a representation $\sigma$ is still denoted by $\sigma$ for simplicity. For more details on representation of groups, we refer the book \cite{Folland}.

Let $A$ be a (complex) Banach space. Let $L^p(G,A),\,1\leq p<\infty$, denote the space of $A$-valued (classes of) functions on $G$ that are strongly Bochner integrable. It is endowed with the norm
$$\|f\|_{L^p}=\left( \int_G \|f(x)\|_A^pdx\right)^{\frac{1}{p}}.$$
Also, let $L^\infty (G,A)$ denote the space of $A$-valued measurable functions on $G$ that are bounded almost everywhere,  endowed with the essential sup-norm:
$$\|f\|_{L^\infty}=\mbox{esssup}(f)=\inf\{C>0 : \|f(x)\|\leq C, \, \mbox{almost everywhere} \}.$$
Let $f\in L^1(G,A)$. The Fourier coefficient of $f$,  denoted by $\widehat{f}(\sigma)$,  is  defined by (see \cite{AssiamouaOlubummo1989}) :
\begin{equation}
\widehat{f}(\sigma)=\int_G f(x)\sigma (x)^*dx, \sigma \in \widehat{G}
\end{equation}
where $\sigma (x)^*$ is the adjoint of the operator $\sigma (x)$ and the integration on $G$ is taken with respect to its normalized Haar measure. 
 Each coefficient  $\widehat{f}(\sigma)$ is interpreted here as the sesquilinear mapping from $H_\sigma \times  H_\sigma$ into $A$ such that 
 \begin{equation}
 \widehat{f}(\sigma)(\xi,\eta)=\int_G \langle \sigma (x)^*\xi,\eta \rangle_{H_\sigma}f(x) dx,\, \xi,\eta\in  H_\sigma.  
\end{equation}  
Consider an orthonormal basis $\{\xi_1^\sigma, \cdots, \xi_{d_\sigma}^\sigma\}$ of $H_{\sigma}$ and define the matrix coefficients  $u_{ij}^\sigma$ by 
$$u_{ij}^\sigma (x)=\langle\sigma (x)\xi_j^\sigma, \xi_i^\sigma \rangle_{H_\sigma},\, x\in G. $$
The inversion formula is given as follows
\begin{equation}
f(x)=\sum\limits_{\sigma \in \widehat{G}}d_\sigma \sum\limits_{i=1}^{d_\sigma}\sum\limits_{j=1}^{d_\sigma}u^\sigma_{ij}(x)\widehat{f}(\sigma)(\xi_j^\sigma, \xi_i^\sigma),\, x\in G.
\end{equation}

Let $\mathscr{S}(H_\sigma \times H_\sigma,A)$ denote the set of sesquilinear maps from $H_\sigma \times H_\sigma$ into $A$. Set  
$$\mathscr{S}(\widehat{G},A)=\prod\limits_{\sigma \in \widehat{G}}\mathscr{S}(H_\sigma \times H_\sigma,A).$$
Now let us recall the definition of  the Assiamoua spaces $\mathscr{S}_p(\widehat{G},A),\, 1\leq p\leq \infty$. For $1\leq p< \infty$,  the space $\mathscr{S}_p(\widehat{G},A)$ is defined  as  the familly of sesquilinear maps $\left(\varphi (\sigma)\right)_{\sigma \in \widehat{G}}\in \mathscr{S}(\widehat{G},A)$ such that 
\begin{equation}
\sum\limits_{\sigma \in \widehat{G}}d_\sigma \sum\limits_{i=1}^{d_\sigma}\sum\limits_{j=1}^{d_\sigma}\|\varphi (\sigma)(\xi_j^\sigma,\xi_i^\sigma)\|_A^p<\infty.
\end{equation}
The space $\mathscr{S}_p(\widehat{G},A)$ is equipped with the norm  defined by
\begin{equation}
\|\varphi\|_{\mathscr{S}_p}=\left(  \sum\limits_{\sigma \in \widehat{G}}d_\sigma \sum\limits_{i=1}^{d_\sigma}\sum\limits_{j=1}^{d_\sigma}\|\varphi (\sigma)(\xi_j^\sigma,\xi_i^\sigma)\|_A^p\right)^{\frac{1}{p}}.
\end{equation}
Define the space  $\mathscr{S}_\infty(\widehat{G},A)$ as the collection ofof sesquilinear maps $\left(\varphi (\sigma)\right)_{\sigma \in \widehat{G}}\in \mathscr{S}(\widehat{G},A)$ such that 
$$
\|\varphi\|_{\mathscr{S}_\infty}:=\sup\{\|\varphi(\sigma)\| : \sigma \in \widehat{G}\}<\infty
$$
where
$$\|\varphi(\sigma)\|=\sup\{\|\varphi(\sigma)(\xi,\eta)\|_A : \|\xi\|_{H_\sigma}\leq 1, \|\eta\|_{H_\sigma}\leq 1\}.$$
 It has been proven that the Assiamoua spaces   are complete spaces \cite{Mensah1}. Particularly, the spaces $\mathscr{S}_p(\widehat{G},A),\, p=1,2$  will be used in this article. The space $\mathscr{S}_1(\widehat{G},A)$ will form the foundation for defining the  spectral Barron spaces which are the subject of this study. 
Finally, to fully grasp  additional properties of the Assiamoua spaces $\mathscr{S}_p(\widehat{G},A)$  and their connection with the Fourier transform described above,  we refer to the  articles \cite{AssiamouaOlubummo1989, Mensah1, Mensah3}.
\section{Functional and embedding results}\label{Functional and embedding results}
\subsection{Definition and examples}
In this section, we list our main results. Mainly, we define the spectral Barron  spaces, study their main functional properties and their connection with the Sobolev spaces associated with this framework.

Let $s$ be a nonnegative real number. Let $\gamma : \widehat{G}\rightarrow (0,\infty)$ be a measurable function. 

\begin{definition}
Let $G$ be a compact group and $A$ a Banach space. 
 We call spectral Barron space (of order $s$) the set 
\begin{equation}
\mathfrak{B}^s_\gamma (G, A) =\left\{ f \in L^2(G,A): (1+\gamma (\sigma)^2)^{\frac{s}{2}}\widehat{f}\in \mathscr{S}_1(\widehat{G}, A)\right\}
\end{equation}
equipped with the norm
\begin{equation}
\|f\|_{\mathfrak{B}^s_\gamma}=\sum_{\sigma\in \widehat{G}}d_\sigma  (1+\gamma (\sigma)^2)^{\frac{s}{2}}\sum_{i=1}^{d_\sigma}\sum_{j=1}^{d_\sigma}\|\widehat{f}(\sigma)(\xi_j^\sigma, \xi_i^\sigma)\|_A=\|(1+\gamma (\sigma)^2)^{\frac{s}{2}}\widehat{f}\|_{\mathscr{S}_1}. 
\end{equation}
\end{definition}

In order to give a concrete idea of spectral Barron spaces in our context, we provide two examples before moving further.
\begin{example}[Spectral Barron spaces on the torus] {\rm 
Let $\mathbb{T}=\{z\in \mathbb{C} : |z|=1\}$ be the one-dimensional torus. It can be  identified in a natural way with the interval $[0, 1)$ of the real line through the group isomorphism $e^{2i\pi x} \rightarrow x$. We know that $\mathbb{T}$ is a compact abelian group.  Thus, all its unitary irreducible representations are  1-dimensional. Its unitary dual $\widehat{\mathbb{T}}$ is identified with the set $\mathbb{Z}$ of integers. Let $A$ denote a Banach space. 
Here, the Fourier coefficients of $f\in L^1(\mathbb{T},A)$ are given by 
\begin{equation}
\widehat{f}(n)=\int_{\mathbb{T}}f(x)e^{-2i\pi nx}dx,\, n\in \mathbb{Z},  
\end{equation}
where the integration is understood in the sense of Bochner.  
Let $\gamma$ be any positive sequence indexed by $\mathbb{Z}$. 
In this context, the spectral Barron space $\mathfrak{B}^s_\gamma (\mathbb{T}, A)$ consists of all functions $f\in L^2(\mathbb{T}, A)$ such that 
$$\|f\|_{\mathfrak{B}^s_\gamma}:=\sum\limits_{n=-\infty}^\infty (1+\gamma (n)^2)^{\frac{s}{2}}\|\widehat{f}(n)\|_A <\infty.$$
}\end{example}

\begin{example}[Spectral Barron spaces on $SU(2)$] {\rm
The linear group $SU(2)$ consists of $2\times 2$-matrices $M$ having coefficients in $\mathbb{C}$ such that $\det (M)=1$ and $M^*M=I$ where $M^*$ is the adjoint of $M$ and $I$ is the identity matrix. It is known that elements  $x$ of $SU(2)$ are   of the form 
\begin{equation}
x=\left(\begin{array}{cc}
\alpha & \beta \\ 
-\overline{\beta} & \overline{\alpha}
\end{array} \right), \alpha,\beta \in \mathbb{C}, \mbox{ with } |\alpha|^2+|\beta|^2=1.
\end{equation}

For $l\in \frac{1}{2}\mathbb{N}_0$, set 
$$V_l=\left\{f\in \mathbb{C}[z_1,z_2] : f(z_1,z_2)=\sum\limits_{k=0}^{2l}a_kz_1^kz_2^{2l-k},\, a_0, \cdots, a_{2l}\in \mathbb{C} \right\}.$$
There is a inner product on $V_l$ defined by 
\begin{equation}
\langle f,g\rangle=\left\langle\sum\limits_{i=0}^{2l}a_iz_1^iz_2^{2l-i}, \sum\limits_{j=0}^{2l}b_jz_1^jz_2^{2l-j} \right\rangle=\sum\limits_{k=0}^{2l}k!(2l-k)!a_k\overline{b_k}.
\end{equation}

Let $U(V_l)$ denote the set of unitary operators on $V_l$. 
The unitary irreducible representations of $SU(2)$ are the mappings $T_l,\, l=0,\frac{1}{2},\frac{2}{2}, \frac{3}{2},  \cdots$ constructed as follows (see \cite{Dooley} or \cite[Chapter 11]{Ruzhansky}) :
$T_l : SU(2)\rightarrow  U(V_l)$ with 

\begin{equation}
\left(T_l\left(\begin{array}{cc}
\alpha & \beta \\ 
-\overline{\beta} & \overline{\alpha}
\end{array} \right)f\right)(z_1,z_2)=f(\alpha z_1-\overline{\beta}z_2,\beta z_1+\overline{\alpha}z_2),\quad f\in V_l.
\end{equation}

Clearly, an  orthonormal basis  of $V_l$ is the collection $\{p_k^l : k=0,\cdots, 2l\}$ where
$p_k^l(z_1,z_2)=z_1^kz_2^{2l-k}$. 
Thus, the  dimension of the vector space $V_l$ is $d_l=2l+1$. Let $A$ be a Banach space. 
The Fourier coefficient $\widehat{f}(l)$ of $f\in L^1(SU(2), A)$ is defined  by 
$$\widehat{f}(l)(p_{j}^l,p_{i}^l)=\int_{SU(2)}\langle T_l(x)^*p_{j}^l,p_{i}^l\rangle f(x) dx.$$
In this context, the spectral Barron space $\mathfrak{B}^s_\gamma(SU(2),A)$ consists of functions $f\in L^2( SU(2),A)$ such that
$$\|f\|_{\mathfrak{B}^s_\gamma}=\sum\limits_{l\in \frac{1}{2}\mathbb{N}_0}(2l+1)(1+\gamma (l)^2)^{\frac{s}{2}}\sum\limits_{i=0}^{2l}\sum\limits_{j=0}^{2l}\|\widehat{f}(l)(p_{j}^l, p_{i}^l)\|_A<\infty.$$ 
}\end{example}

\subsection{Functional properties of spectral Barron spaces}
The spectral Barron spaces $\mathfrak{B}^s_\gamma (G, A)$ may not be complete. However, as a normed space, it admits a unique continuous completion which we denote by  $\tilde{\mathfrak{B}}^s_\gamma (G, A)$.

Now, we consider the operator $(I-\Delta)^s$ defined by 
\begin{eqnarray}
(I-\Delta)^sf(x)&=\sum\limits_{\sigma\in \widehat{G}}d_\sigma  (1+\gamma(\sigma)^2)^s\sum\limits_{i=1}^{d_\sigma}\sum\limits_{j=1}^{d_\sigma}\widehat{f}(\sigma)(\xi_j^\sigma, \xi_i^\sigma)u_{i,j}^\sigma (x)\\
&=\mathscr{F}^{-1}((1+\gamma(\sigma)^2)^s\widehat{f}
), \quad f\in  \mathfrak{B}^{2s}_\gamma (G, A).
\end{eqnarray}

\begin{theorem}\label{isomorphism}
Let $s>0$. The operator $(I-\Delta)^s$ is an isometric isomorphism from $\mathfrak{B}^{2s}_\gamma (G, A)$ onto $\mathfrak{B}^{0}_\gamma (G, A)$.
\end{theorem}
\begin{proof}
Clearly, for $f\in \mathfrak{B}^{2s}_\gamma (G, A)$, we have 
$\|(I-\Delta)^sf\|_{\mathfrak{B}^0_\gamma}=\|f\|_{\mathfrak{B}^{2s}_\gamma}$  and the injectivity of the operator $(I-\Delta)^s$ follows. 
Now, let $g$ be in $\mathfrak{B}^0_\gamma (G,A)$. Then, $g\in L^2(G,A)$ and $\widehat{g}\in \mathscr{S}_1(\widehat{G}, A)$.  In fact, $\widehat{g}\in \mathscr{S}_2(\widehat{G}, A)$

 Set $h=\mathscr{F}^{-1}\left(\displaystyle\frac{\widehat{g}}{(1+\gamma(\sigma)^2)^s}\right)$. Then, $h\in \mathfrak{B}^{2s}_\gamma (G, A)$ and  $(1-\Delta)^sh=g$. Thus, the operator $(I-\Delta)^s$ is surjective. 
\end{proof}

The following theorem is an interpolation  result between spectral Barron spaces. 
\begin{theorem}
Let $0\leq r\leq t$. If $f\in \mathfrak{B}^r_\gamma (G, A)$, then $\forall s\in [r,t]$ with $s=\alpha r+(1-\alpha)t, \, 0\leq \alpha \leq 1$, we have 
\begin{equation}
\|f\|_{\mathfrak{B}^s_\gamma}\leq \|f\|_{\mathfrak{B}^r_\gamma}^\alpha\|f\|_{\mathfrak{B}^t_\gamma}^{1-\alpha}.
\end{equation}
\end{theorem}
\begin{proof}
\begin{align*}
\|f\|_{\mathfrak{B}^s_\gamma}&= \sum_{\sigma\in \widehat{G}}d_\sigma  (1+\gamma (\sigma)^2)^{\frac{s}{2}}\sum_{i=1}^{d_\sigma}\sum_{j=1}^{d_\sigma}\|\widehat{f}(\sigma)(\xi_j^\sigma, \xi_i^\sigma)\|_A\\
&=\sum_{\sigma\in \widehat{G}}d_\sigma  (1+\gamma (\sigma)^2)^{\frac{\alpha r}{2}}(1+\gamma (\sigma)^2)^{\frac{(1-\alpha)t}{2}}\sum_{i=1}^{d_\sigma}\sum_{j=1}^{d_\sigma}\|\widehat{f}(\sigma)(\xi_j^\sigma, \xi_i^\sigma)\|_A^\alpha\|\widehat{f}(\sigma)(\xi_j^\sigma, \xi_i^\sigma)\|_A^{1-\alpha}\\
&\leq \left( \sum_{\sigma\in \widehat{G}}d_\sigma  (1+\gamma (\sigma)^2)^{\frac{ r}{2}}  \sum_{i=1}^{d_\sigma}\sum_{j=1}^{d_\sigma}\|\widehat{f}(\sigma)(\xi_j^\sigma, \xi_i^\sigma)\|_A\right)^{\alpha} \\
&\times\left( \sum_{\sigma\in \widehat{G}}d_\sigma  (1+\gamma (\sigma)^2)^{\frac{ t}{2}}  \sum_{i=1}^{d_\sigma}\sum_{j=1}^{d_\sigma}\|\widehat{f}(\sigma)(\xi_j^\sigma, \xi_i^\sigma)\|_A\right)^{1-\alpha}\\
&\leq \|f\|_{\mathfrak{B}^r_\gamma}^\alpha \|f\|_{\mathfrak{B}^t_\gamma}^{1-\alpha}.
\end{align*}
\end{proof}

Let $s\geq 0$ and $t\geq 0$. Consider a measurable function $a: \widehat{G}\rightarrow A$.  Consider the pseudo-differential operator $\mathfrak{P}$ with symbol $a$ defined by
\begin{equation}
\mathfrak{P}f(x)=\sum_{\sigma\in \widehat{G}}d_\sigma  a(\sigma)\sum_{i=1}^{d_\sigma}\sum_{j=1}^{d_\sigma}\widehat{f}(\sigma)(\xi_j^\sigma, \xi_i^\sigma)u_{i,j}^\sigma (x).
\end{equation}
\begin{theorem}
If $(1+\gamma (\sigma)^2)^{\frac{s-t}{2}}a \in L^{\infty}(\widehat{G}, A)$, then $\mathfrak{P} : \mathfrak{B}^{t}_\gamma (G, A) \rightarrow \mathfrak{B}^{s}_\gamma (G, A)$ is a bounded linear operator.
\end{theorem}
\begin{proof}
\begin{align*}
\|\mathfrak{P}f\|_{\mathfrak{B}^{s}_\gamma}&=\sum_{\sigma\in \widehat{G}}d_\sigma  (1+\gamma (\sigma)^2)^{\frac{s}{2}}\sum_{i=1}^{d_\sigma}\sum_{j=1}^{d_\sigma}\|\widehat{\mathfrak{P}f}(\sigma)(\xi_j^\sigma, \xi_i^\sigma)\|_A\\
&=\sum_{\sigma\in \widehat{G}}d_\sigma  (1+\gamma (\sigma)^2)^{\frac{s}{2}}|a(\sigma)|\sum_{i=1}^{d_\sigma}\sum_{j=1}^{d_\sigma}\|\widehat{f}(\sigma)(\xi_j^\sigma, \xi_i^\sigma)\|_A\\
&=\sum_{\sigma\in \widehat{G}}d_\sigma  (1+\gamma (\sigma)^2)^{\frac{s-t}{2}}|a(\sigma)|(1+\gamma (\sigma)^2)^{\frac{t}{2}}\sum_{i=1}^{d_\sigma}\sum_{j=1}^{d_\sigma}\|\widehat{f}(\sigma)(\xi_j^\sigma, \xi_i^\sigma)\|_A\\
&\leq \|(1+\gamma (\cdot)^2)^{\frac{s-t}{2}}a\|_{L^\infty}\sum_{\sigma\in \widehat{G}}d_\sigma  (1+\gamma (\sigma)^2)^{\frac{t}{2}}\sum_{i=1}^{d_\sigma}\sum_{j=1}^{d_\sigma}\|\widehat{f}(\sigma)(\xi_j^\sigma, \xi_i^\sigma)\|_A\\
&=\|(1+\gamma (\cdot)^2)^{\frac{s-t}{2}}a\|_{L^\infty}\|f\|_{\mathfrak{B}^{t}_\gamma}.
\end{align*}
The conclusion follows.
\end{proof}

In the following result, we need $A$ to be a Banach algebra in order to have the possibility to multiply  its elements. We recall that a Banach algebra is a Banach space $A$ with an internal associative product which is compatible with vector space structure of $A$ such that for all $x,y\in A$, 
$\|xy\|\leq \|x\|\|y\|$.

We assume, only in this part, that $A$ is a Banach algebra.  For $f,g \in L^1(G,A)$,  the convolution product $f\ast g$ is defined by 
\begin{equation}
f\ast g (x)=\int_G f(y)g(y^{-1}x)dy,\quad x\in G.
\end{equation}
Then, it  follows that \cite{Assiamoua2}:
\begin{equation}
\widehat{f\ast g}(\sigma)(\xi_j^\sigma,\xi_i^\sigma)=\sum\limits_{k=1}^{d\sigma}\widehat{f}(\sigma)(\xi_k^\sigma,\xi_i^\sigma)\widehat{g}(\sigma)(\xi_j^\sigma,\xi_k^\sigma). 
\end{equation}
\begin{theorem}
Assume that $A$ is a Banach algebra and  that $\rho:=\sup\{d_\sigma :\sigma\in \widehat{G}\}<\infty$.  If $f\in L^1(G,A)$ and  $g\in \mathfrak{B}^{s}_\gamma(G,A)$, then $f\ast g\in   \mathfrak{B}^{s}_\gamma(G,A)$ and 
$$\|f\ast g\|_{\mathfrak{B}^{s}_\gamma}\leq \rho \|f\|_{L^1}\|g\|_{\mathfrak{B}^{s}_\gamma}.$$
\end{theorem}
\begin{proof}
\begin{align*}
\left\|\widehat{f\ast g}(\sigma)(\xi_j^\sigma,\xi_i^\sigma)\right\|_A &=\left\|\sum\limits_{k=1}^{d_\sigma}\widehat{f}(\sigma)(\xi_k^\sigma,\xi_i^\sigma)\widehat{g}(\sigma)(\xi_j^\sigma,\xi_k^\sigma)\right\|_A\\
&\leq  \sum\limits_{k=1}^{d_\sigma}\left\|\widehat{f}(\sigma)(\xi_k^\sigma,\xi_i^\sigma)\right\|_A\left\|\widehat{g}(\sigma)(\xi_j^\sigma,\xi_k^\sigma)\right\|_A. 
\end{align*}
Set  $\Gamma=\sum\limits_{\sigma\in \widehat{G}}d_\sigma (1+\gamma (\sigma)^2)^{\frac{s}{2}}\sum\limits_{i=1}^{d_\sigma}\sum\limits_{j=1}^{d_\sigma}\left\|\widehat{f\ast g}(\sigma)(\xi_j^\sigma,\xi_i^\sigma)\right\|_A$.
This implies, 
\begin{align*}
\Gamma&\leq \sum_{\sigma\in \widehat{G}}d_\sigma (1+\gamma (\sigma)^2)^{\frac{s}{2}} \sum_{i=1}^{d_\sigma}\sum_{j=1}^{d_\sigma}\sum\limits_{k=1}^{d_\sigma}
\left\|\widehat{f}(\sigma)(\xi_k^\sigma,\xi_i^\sigma)\right\|_A
\left\|\widehat{g}(\sigma)(\xi_j^\sigma,\xi_k^\sigma)\right\|_A\\
&\leq \|\widehat{f}\|_{\mathscr{S}_\infty}\sum_{\sigma\in \widehat{G}}d_\sigma^2 (1+\gamma (\sigma)^2)^{\frac{s}{2}} \sum_{j=1}^{d_\sigma}\sum\limits_{k=1}^{d_\sigma}
\left\|\widehat{g}(\sigma)(\xi_j^\sigma,\xi_k^\sigma)\right\|_A\\
&\leq \rho\|\widehat{f}\|_{\mathscr{S}_\infty}\sum_{\sigma\in \widehat{G}}d_\sigma (1+\gamma (\sigma)^2)^{\frac{s}{2}} \sum_{j=1}^{d_\sigma}\sum\limits_{k=1}^{d_\sigma}
\left\|\widehat{g}(\sigma)(\xi_j^\sigma,\xi_k^\sigma)\right\|_A\\
&=\rho\|\widehat{f}\|_{\mathscr{S}_\infty}\|g\|_{\mathfrak{B}^{s}_\gamma}\\
&\leq \rho\|f\|_{L^1}\|g\|_{\mathfrak{B}^{s}_\gamma}<\infty.
\end{align*}
That is, $f\ast g$ belongs to $\mathfrak{B}^{s}_\gamma (G,A)$ and  $\|f\ast g\|_{\mathfrak{B}^{s}_\gamma}\leq \rho \|f\|_{L^1}\|g\|_{\mathfrak{B}^{s}_\gamma}.$
\end{proof}

\subsection{Continuous embedding results for spectral Barron spaces}
The spectral Barron spaces  embed continuously in other function spaces. 
The first result states continuous embedding between spectral Barron spaces. 
\begin{theorem}
Let $0\leq s<t$. Then, $\mathfrak{B}^{t}_\gamma (G,A)$ embeds continuously into  $\mathfrak{B}^{s}_\gamma (G,A)$ and $$\forall f\in \mathfrak{B}^{s}_\gamma (G,A),\, 
\|f\|_{\mathfrak{B}^{s}_\gamma} \leq \|f\|_{\mathfrak{B}^{t}_\gamma}.$$ 
\end{theorem}
\begin{proof}
The proof is trivial by noting that since  $0\leq s<t$ and $1+\gamma (\sigma)^2>1$, then 
$$(1+\gamma (\sigma)^2)^{\frac{s}{2}}< (1+\gamma (\sigma)^2)^{\frac{t}{2}}.$$

\end{proof}

 Now, we will connect the spectral Barron space with the Sobolev space  $H^s_\gamma (G,A)$ of vector-valued functions on  a compact group  studied in the article \cite{Mensah2}. Let us recall that this Sobolev space is defined as follows:

\begin{equation}
H^s_\gamma (G,A)=\left\{f\in L^2(G,A) : (1+\gamma (\sigma)^2)^{\frac{s}{2}}\widehat{f}\in \mathscr{S}_2(\widehat{G}, A)\right\}.
\end{equation}
In other words,  $f\in  L^2(G,A)$ belongs to  $H^s_\gamma (G,A)$ if $$\sum_{\sigma\in \widehat{G}}d_\sigma  (1+\gamma (\sigma)^2)^s\sum_{i=1}^{d_\sigma}\sum_{j=1}^{d_\sigma}\|\widehat{f}(\sigma)(\xi_j^\sigma, \xi_i^\sigma)\|_A^2<\infty.$$

Naturally, the norm on $H^s_\gamma (G,A)$ is given by

\begin{equation}
\|f\|_{H^s_\gamma}=\left(  \sum_{\sigma\in \widehat{G}}d_\sigma  (1+\gamma (\sigma)^2)^s\sum_{i=1}^{d_\sigma}\sum_{j=1}^{d_\sigma}\|\widehat{f}(\sigma)(\xi_j^\sigma, \xi_i^\sigma)\|_A^2 \right)^{\frac{1}{2}}. 
\end{equation}
The following theorem provides a sufficient condition for the continuous embedding of  Sobolev spaces into spectral Barron spaces. 
\begin{theorem} Let $0\leq s< t$. 
If $\kappa:=\sum\limits_{\sigma\in \widehat{G}}d_\sigma  (1+\gamma (\sigma)^2)^{t-s}<\infty$, then $H^s_\gamma (G,A)$ embeds continuously into $\mathfrak{B}^{t}_\gamma(G,A)$ and $$\forall f\in H^s_\gamma (G,A), \, \|f\|_{\mathfrak{B}^{t}_\gamma}\leq \kappa^{\frac{1}{2}}\|f\|_{H^{s}_\gamma}.$$
\end{theorem}
\begin{proof}
\begin{align*}
\|f\|_{\mathfrak{B}^{t}_\gamma}&=\sum_{\sigma\in \widehat{G}}d_\sigma  (1+\gamma (\sigma)^2)^{\frac{t}{2}}\sum_{i=1}^{d_\sigma}\sum_{j=1}^{d_\sigma}\|\widehat{f}(\sigma)(\xi_j^\sigma, \xi_i^\sigma)\|_A\\
&=\sum_{\sigma\in \widehat{G}}d_\sigma  (1+\gamma (\sigma)^2)^{\frac{t-s}{2}}(1+\gamma (\sigma)^2)^{\frac{s}{2}}\sum_{i=1}^{d_\sigma}\sum_{j=1}^{d_\sigma}\|\widehat{f}(\sigma)(\xi_j^\sigma, \xi_i^\sigma)\|_A\\
&\leq  \left(\sum_{\sigma\in \widehat{G}}d_\sigma  (1+\gamma (\sigma)^2)^{t-s}\right)^{\frac{1}{2}}\left(   \sum_{\sigma\in \widehat{G}}d_\sigma (1+\gamma (\sigma)^2)^{s}\sum_{i=1}^{d_\sigma}\sum_{j=1}^{d_\sigma}\|\widehat{f}(\sigma)(\xi_j^\sigma, \xi_i^\sigma)\|_A^2 \right)^{\frac{1}{2}}\\
&(\mbox{by the H\"older inequality})\\
&=\kappa^{\frac{1}{2}}\|f\|_{H^{s}_\gamma}.
\end{align*}
\end{proof}

The following theorem states the continuous embedding of $\mathfrak{B}^{s}_\gamma (G,A)$ into  $L^{\infty}(G,A)$.
\begin{theorem}
Let $s\geq0$. Then, $\mathfrak{B}^{s}_\gamma (G,A)$ embeds continuously into $L^{\infty}(G,A)$ and $$\forall f\in \mathfrak{B}^{s}_\gamma (G,A), \, \|f\|_{L^\infty}\leq \|f\|_{\mathfrak{B}^{s}_\gamma}.$$
\end{theorem}
\begin{proof}
\begin{align*}
\|f(x)\|_A&= \left\|\sum\limits_{\sigma \in \widehat{G}}d_\sigma \sum\limits_{i=1}^{d_\sigma}\sum\limits_{j=1}^{d_\sigma}\widehat{f}(\sigma)(\xi_j^\sigma,\xi_i^\sigma)u_{ij}^\sigma (x)\right\|_A\\
&\leq \sum\limits_{\sigma \in \widehat{G}}d_\sigma \sum\limits_{i=1}^{d_\sigma}\sum\limits_{j=1}^{d_\sigma}\left\|\widehat{f}(\sigma)(\xi_j^\sigma,\xi_i^\sigma)\right\|_A |u_{ij}^\sigma (x)|\\
&\leq \sum\limits_{\sigma \in \widehat{G}}d_\sigma \sum\limits_{i=1}^{d_\sigma}\sum\limits_{j=1}^{d_\sigma}\left\|\widehat{f}(\sigma)(\xi_j^\sigma,\xi_i^\sigma)\right\|_A\\
&\leq \sum\limits_{\sigma \in \widehat{G}}d_\sigma (1+\gamma (\sigma)^2)^{\frac{s}{2}}\sum\limits_{i=1}^{d_\sigma}\sum\limits_{j=1}^{d_\sigma}\left\|\widehat{f}(\sigma)(\xi_j^\sigma,\xi_i^\sigma)\right\|_A\\
&=\|f\|_{\mathfrak{B}^{s}_\gamma}.
\end{align*}
Thus, $\|f\|_{L^\infty}\leq \|f\|_{\mathfrak{B}^{s}_\gamma}$. 
\end{proof}

\section{A Barron-type approximation theorem}\label{approximation}
 We assume in this section that $A$ has a  Hilbert space structure in order to use Maurey's empirical method \cite{Jones1992, Pisier1981}.

For $f \in \mathfrak{B}^s_\gamma(G,A)$,  write
$$
M := \|f\|_{\mathfrak{B}^s_\gamma}
= \sum_{\sigma \in \widehat{G}} d_\sigma (1+\gamma(\sigma)^2)^{s/2}
  \sum_{i=1}^{d_\sigma}\sum_{j=1}^{d_\sigma}
  \big\| \widehat{f}(\sigma)(\xi^\sigma_j,\xi^\sigma_i) \big\|_A < \infty .
$$
 
We use the elementary fact that matrix coefficients of unitary
representations  are uniformly bounded by $1$. Indeed, 
$u^\sigma_{ij}(x) = \langle \sigma(x)\xi^\sigma_j,\xi^\sigma_i\rangle_{H_\sigma}$ with
$\sigma(x)$ unitary and $\xi^\sigma_i,\xi^\sigma_j$ unit vectors in $H_\sigma$. The  Cauchy-Schwarz inequality gives
\begin{equation}\label{eq:bounded-atoms}
|u^\sigma_{ij}(x)| \le \|\sigma(x)\xi^\sigma_j\|_{H_\sigma}\,\|\xi^\sigma_i\|_{H_\sigma} = 1,
\quad \forall x \in G,\ \forall \sigma \in \widehat{G}.
\end{equation}
 
For every index $k=(\sigma,i,j)\in \mathfrak{K}:=\widehat{G}\times \{1,\cdots,d_\sigma\}^2$ with $\widehat{f}(\sigma)(\xi^\sigma_j,\xi^\sigma_i)\ne 0$, set
$$
c_k := \big\|\widehat{f}(\sigma)(\xi^\sigma_j,\xi^\sigma_i)\big\|_A, \qquad
e_k := \frac{\widehat{f}(\sigma)(\xi^\sigma_j,\xi^\sigma_i)}{c_k} \in A.$$ Then, $\|e_k\|_A=1$.
 Set 
 $$
a_k := d_\sigma (1+\gamma(\sigma)^2)^{s/2} c_k \ge 0, \qquad
\psi_k(x) := (1+\gamma(\sigma)^2)^{-s/2}\, u^\sigma_{ij}(x)\, e_k.
$$
By definition $\sum\limits_{k\in\mathfrak{K}} a_k = M$, and by \eqref{eq:bounded-atoms} together with $s \ge 0$, we have 
$$
\|\psi_k(x)\|_A = (1+\gamma(\sigma)^2)^{-s/2}\,|u^\sigma_{ij}(x)| \le 1
\quad \text{for all } x \in G.
$$
Therefore,  $\|\psi_k\|_{L^\infty(G,A)} \le 1$. Hence $\|\psi_k\|_{L^2(G,A)}\le 1$.
 
The inversion formula
$$
f(x) = \sum_{\sigma\in\widehat{G}} d_\sigma \sum_{i=1}^{d_\sigma}\sum_{j=1}^{d_\sigma}
u^\sigma_{ij}(x)\widehat{f}(\sigma)(\xi^\sigma_j,\xi^\sigma_i)
$$
now reads
\begin{equation}\label{eq:barycenter}
f = \sum_{k\in \mathfrak{K}} a_k\, \psi_k. 
\end{equation} 
Since $\sum\limits_{k\in\mathfrak{K}} \frac{a_k}{M} = 1$, 
the function  $\frac{f}{M}$ is a convex combination  of the bounded
set $\{\psi_k\} \subset L^2(G,A)$.

\begin{theorem}\label{thm:main}
Let $G$ be a compact group and $A$ a separable complex Hilbert space. Let  $f \in \mathfrak{B}^s_\gamma(G,A)$ and $M = \|f\|_{\mathfrak{B}^s_\gamma}$. For every $n \in \mathbb{N}$, 
there exist indices $(\sigma_\ell,i_\ell,j_\ell)_{\ell=1}^n$ (with possible
repetitions) and unit vectors $e_1,\dots,e_n \in A$ such that the $n$-term
matrix-coefficient approximant
\[
f_n(x) := \frac{M}{n} \sum_{\ell=1}^n
(1+\gamma(\sigma_\ell)^2)^{-s/2}\, u^{\sigma_\ell}_{i_\ell j_\ell}(x)\, e_\ell
\]
satisfies
$$
\|f - f_n\|_{L^2(G,A)} \;\le\; \frac{M}{\sqrt n}.
$$
\end{theorem}
 
\begin{proof}
Set $p_k = \frac{a_k}{M} \ge 0$, so $\sum\limits_{k\in\mathfrak{K}} p_k = 1$.
Let $K_1,\dots,K_n$ be i.i.d.\ random indices with probability density $\mathbb{P}(K_\ell=k)=p_k$,
and  $Y_\ell = \psi_{K_\ell} \in L^2(G,A)$, i.i.d.
 
By construction $\mathbb{E}[Y_1] = \sum\limits_{k\in\mathfrak{K}} p_k \psi_k = \frac{f}{M}$. Since
$Y_1,\dots,Y_n$ are i.i.d.\ elements of the Hilbert space $L^2(G,A)$, we have
\begin{align*}
\mathbb{E}\left\| \frac{1}{n}\sum_{\ell=1}^n Y_\ell - \frac{f}{M} \right\|_{L^2(G,A)}^2
&= \frac1n \Var(Y_1)\\
&= \frac1n \Big( \mathbb{E}\|Y_1\|_{L^2(G,A)}^2 - \|\mathbb{E} Y_1\|_{L^2(G,A)}^2 \Big)\\
&\le \frac1n\, \mathbb{E}\|Y_1\|_{L^2(G,A)}^2.
\end{align*}
Since $\|Y_1\|_{L^2(G,A)} = \|\psi_{K_1}\|_{L^2(G,A)} \le 1$ almost surely,
\[
\mathbb{E}\left\| \frac1n\sum_{\ell=1}^n Y_\ell - \frac{f}{M} \right\|_{L^2(G,A)}^2 \le \frac1n .
\]
As this holds in expectation, some realization $(k_1,\dots,k_n)$ of
$(K_1,\dots,K_n)$ attains at most the average value, i.e.
\[
\left\| \frac1n \sum_{\ell=1}^n \psi_{k_\ell} - \frac{f}{M} \right\|_{L^2(G,A)}^2
\le \frac1n .
\]
Multiplying by $M$ and setting $f_n := \dfrac{M}{n}\sum_{\ell=1}^n \psi_{k_\ell}$
gives $\|f-f_n\|_{L^2(G,A)} \le\frac{M}{\sqrt {n}}$. 
\end{proof}
 \begin{remark}{\rm
 \begin{enumerate}
\item[(i)] {\bf Dimension independence.}
The constant $\frac{M}{\sqrt {n}}$ depends only on $\|f\|_{\mathfrak{B}^s_\gamma}$ (not on $G$,
not on $\dim A$, and not on the growth of $d_\sigma$ or of the weight
$\gamma$. This is the direct analogue of the dimension-independent rate in
Barron's classical theorem on $\mathbb{R}^d$ \cite{Barron1993}, with the matrix coefficients
$u^\sigma_{ij}$ playing the role of the bounded ridge functions
$\sigma(w\cdot x + b)$.
\item[(ii)] {\bf Extension beyond Hilbert space.}
The Hilbert space structure of $A$ is used only through the variance
identity $\mathbb{E}\left\|\frac{1}{n}\sum_{\ell=1}^n Y_\ell - \frac{f}{M} \right\|_{L^2(G,A)}^2
= \frac{1}{n} \Var(Y_1)$. If $B$ is a
Banach space of (Rademacher) type $2$ with type-$2$ constant $T_2(B)$, the
same argument, using Pisier's version of Maurey's empirical method \cite{Pisier1989},  yields
$$
\|f - f_n\|_{L^2(G,B)} \le T_2(B)\, \frac{M}{\sqrt n}.
$$
One recovers Theorem~\ref{thm:main} when $B=A$ is Hilbert, since Hilbert
spaces have $T_2(A)=1$.
\end{enumerate}
}\end{remark}

%%%% Bibliography  %%%%%%%%%%

\end{document}